\documentclass[12 pt]{amsart}

\usepackage{graphicx}
\usepackage{amssymb,latexsym}
\usepackage{amsmath,amsfonts,amstext}

\newtheorem{theorem}{Theorem}
\newtheorem{lemma}{Lemma}

\begin{document}
\title[A FOURIER RESTRICTION ESTIMATE   IN $\mathbb{R}^6$]
      {A Fourier  restriction estimate for surfaces of positive curvature in $\mathbb{R}^6$}
      
\author{Faruk Temur}
\address{Department of Mathematics\\
        University of Illinois at Urbana-Champaign
        Urbana, IL 61820}
\email{temur1@illinois.edu}
\thanks{The author was supported by NSF grant DMS-0900865. The author thanks to his advisor M.B. Erdo\~{g}an for the financial support}
\date{August 5, 2011}    

\begin{abstract}
We improve the best known exponent for the restriction conjecture in $\mathbb{R}^6$  improving  the recent results of Bourgain and Guth.  The proof is  applicable to any dimension n satisfying $n \equiv 0 \mod 3 $ 
\end{abstract}  

\maketitle

\section{Introduction}\label{intro}

In \cite{st} Stein posed the following well-known conjecture. Let $S \subset \mathbb{R}^n$ be a smooth, compact hypersurface with positive definite second fundamental form and $\sigma$ its surface measure. Then for $p>2n/n-1$ and $f \in L^{\infty}(S,\sigma)$ we have
\begin{equation}\label{eq1.1}   
\|\widehat{fd\sigma} \|_p \leq C_{p,S}\|f \|_{\infty}.
\end{equation}

This conjecture is related to some other important problems in harmonic analysis and PDE such as the Kakeya conjecture, the Bochner-Riesz conjecture, local smoothing problem; see \cite{b1, la, t1, w1}. For $n=2$ it is known to be true; see \cite{f}. For $n \geq 3$ it is open despite much effort. The first progress towards this case was the Tomas-Stein theorem \cite{tomas}, and gives $p >2(n+1)/(n-1)$; see \cite{tomas}. Then in  \cite{b}  Bourgain  was able to go below this exponent. Wolff improved Bourgain's result to $(2n^2+n+6)/(n^2+n-1)$, see \cite{w}. Then in three dimensions Tao, Vargas, and Vega further lowered this exponent, and more importantly they developed the bilinear approach which related this conjecture to restriction estimates for compact, transverse subsets of hypersurfaces; see \cite{tvv, tv1}. The work of Tao in \cite{t}, which was a bilinear estimate for compact transverse subsets of paraboloids, through this bilinear method,  verified the conjecture for $p > 2(n+2)/n$. This exponent is the  best one can 
obtain from that approach.

In \cite{bct}, Bennett, Carbery and Tao posed a multilinear version of the restriction conjecture and resolved it. Let $S_1,\ldots,S_m \subset \mathbb{R}^n$ be  smooth compact hypersurfaces that are transverse, that is for any choice of points $\{x_i \in S_i\}$ we have $|x_1' \wedge \cdots \wedge x_m'| >c$ where $x_i'$ is the unit normal at $x_i$, and $c$ some positive constant. Let $\sigma_i$ be the surface measure of $S_i$.  Then for $q >2m/m-1$ and $p' \leq q(m-1)/m$  the result of  \cite{bct} implies
\begin{equation}\label{eq1.2}
\|\prod_{i=1}^m \widehat{f_id\sigma_i}\|_{L^{q/m}(\mathbb{R}^n)} \lesssim \prod_{i=1}^m \|f_i\|_{L^p(S_i)}. 
\end{equation} 
This result by itself does not imply any progress towards the restriction conjecture. But recently in \cite{bg}, Bourgain and Guth combining this with the idea of rescaling, significantly improved the known exponents for $n>4$. The exponents given in \cite{bg}  are  for  every $n>2$ as follows:
\begin{equation}\label{eq1.4}
\begin{aligned}
p>\frac{8n+6}{4n-3} \hspace{3mm} \text{if} \hspace{3mm}n &\equiv 0 \mod3\\
p>\frac{2n+1}{n-1} \hspace{3mm}\text{if}\hspace{3mm} n &\equiv 1 \mod3\\
p>\frac{4n+4}{2n-1} \hspace{3mm}\text{if}\hspace{3mm} n &\equiv 2 \mod 3
\end{aligned} 
\end{equation}
For $n=3,4$ this does not give any improvement. But refining their analysis and combining with Wolff's Kakeya maximal function estimate in \cite{w}, Bourgain and Guth, for $n=3$,  improved the known exponent $p>10/3$ to  $p>33/10$. The aim of this paper is to show that this refined method can be used for $n=6$ too, and to calculate the improvement. We shall also make it clear how to use this strategy for any dimension $n \equiv 0 \mod 3$, though as the improvement is very small and process very technical we will not calculate the improvement for general $n$. We state the case in $n=6$  as a theorem:
\begin{theorem}\label{t1.1}
 Let $S \subset \mathbb{R}^6$ be a smooth, compact hypersurface with positive definite second fundamental form and $\sigma$ its surface measure. Then for $p> 18/7-2/735$ and $f \in L^{\infty}(S,\sigma)$ we have
\begin{equation}\label{eq1.3}   
\|\widehat{fd\sigma} \|_p \lesssim \|f \|_{\infty}.
\end{equation}
\end{theorem} 
Thus the improvement we have over Bourgain-Guth exponent is $2/735$.

The rest of the paper is organized as follows. In section 2, we describe the proof by Bourgain and Guth of (\ref{eq1.4}),  and point out what allows us  when $n\equiv 0 \mod3$ to improve this. In the next section we calculate explicitly the improvement for $n=6$, and at the end of that section, it will be clear to the reader that the process can be repeated to obtain improvement for any $n$ with $n\equiv 0 \mod3$.

\section{The Bourgain-Guth argument}

We first remark that standard $\epsilon$-removal arguments allow us to derive the Theorem 1 from the following theorem:

\begin{theorem}
 Let $S \subset \mathbb{R}^6$ be a smooth, compact hypersurface with positive definite second fundamental form and $\sigma$ its surface measure. Then for $p \geq 18/7-2/735$, and $f \in L^{\infty}(S,\sigma)$ we have
\begin{equation}\label{eq1.3'}   
\|\widehat{fd\sigma} \|_{L^p(B(0,R))} \lesssim R^{\epsilon}\|f \|_{\infty}.
\end{equation}
\end{theorem}

While working in this localized setting we shall use the following version of \eqref{eq1.2} proved in \cite{bct}: for $S_1,S_2,\ldots,S_m$ satisfying the same properties as described in section 1,  one has for every $\epsilon>0$

\begin{equation}\label{eq1.6'}
\|\prod_{i=1}^{m} \widehat{f_id\sigma_i}\|_{L^{q/m}(B(0,R))} \lesssim R^{\epsilon}\prod_{i=1}^{m}\|f_i\|_{L^2(S_i)}
\end{equation}
for $q \geq 2m/m-1$.

We further remark that it suffices to prove  Theorem 2 for $\|f\|_{\infty} \leq 1$, and we accordingly define $Q_R^p$ to be the best constant satisfying
\begin{equation}   
\|\widehat{fd\sigma} \|_{L^p(B(0,R))} \leq Q^p_R.
\end{equation}
This constant clearly is well defined by the crude estimate 

\[\|\widehat{fd\sigma} \|_{L^p(B(0,R))} \lesssim R^{n/p}.\]
 Thus we reduce to showing that $Q_R^p \lesssim R^{\epsilon}$.

 We continue with several lemmas, proofs of which can be found in \cite{bg}. The first two  lemmas rely on multilinear estimates of \cite{bct}, while the third uses rescaling.

Let S be a compact, smooth hypersurface in $\mathbb{R}^n$  with positive definite second fundamental form. Let for $x \in S$, $x'\in S^{n-1}$ denote the unit normal to the surface at the point $x$, and let $\Gamma: S^{n-1} \rightarrow S$ be the Gauss map. Hence $\Gamma(x')=x$. In what follows we will use the notation $\oint_E$ to denote the average over the set $E$. Now we are ready to state our first lemma.  
\begin{lemma}\label{l2.1} 
Let $U_i \subset S, 1\leq i \leq n $ be small caps such that $|x_1' \wedge \ldots \wedge x_n'| > c $ for all $x_i \in U_i$. Let $D_i \subset U_i$, $1 \leq i\leq n$ be discrete sets of $1/M$-separated points for $M$ large. Then for a bounded function $a$ on $S$
\begin{equation*}
\oint_{B_M} \prod_{i=1}^n|\sum_{\xi \in D_i}a(\xi)e^{-ix \cdot \xi}|^{2/n-1} \ll M^{\epsilon}\prod_{i=1}^{n}(\sum_{\xi \in D_i}|a(\xi)|^2)^{1/{n-1}}  
 \end{equation*}
 where $B_M \subset \mathbb{R}^n$ is a ball of radius $M$.
 \end{lemma}
This lemma is a discretized version of \eqref{eq1.2}, using uncertainty principle one replaces discrete sums with integrals of functions that are constant on $1/M$ neighborhoods of points $\xi$, then applies \eqref{eq1.2}. For details see \cite{bg}.

\begin{lemma}\label{l2.2}
Let $2\leq m \leq n$ and $V$ be a subspace of $\mathbb{R}^n$  of dimension $m$. Let $P_1 \ldots P_m \in S$ be points that satisfy $P_i' \in V$ for all $1 \leq i \leq m$ and $|P_1' \wedge \ldots \wedge P_m'| > c$. Let $U_1, \ldots, U_m \subset S$ be small neighborhoods of $P_1, \ldots , P_m$. Let $M$ be large and $D_i \subset U_i$ be sets of $1/M$ separated points $\xi $ that obey the condition $\text{dist}(\xi',V)< c/M$. Then for $f_i \in L^{\infty}(U_i)$, we have
\[ \oint_{B_M}\prod_{i=1}^m|\sum_{\xi \in D_i}\int_{|\zeta-\xi|<\frac{c}{M}} f_i(\zeta)e^{-ix \cdot \zeta}d\sigma(\zeta)|^{\frac{2}{m-1}}dx \ll \] \[ M^{\epsilon}\big(\oint_{B_M}\prod_{i=1}^m(\sum_{\xi \in D_i}|\int_{|\zeta-\xi|<\frac{c}{M}} f_i(\zeta)e^{-ix \cdot \zeta}d\sigma(\zeta)|^2)^{\frac{1}{2m}}dx \big)^{\frac{2m}{m-1}}. \]
\end{lemma}
To prove this lemma we first use uncertainty principle to discretize it as in Lemma 1, then by the hypothesis $dist(\xi',V)<c/M$ and uncertainty principle do a dimension reduction to $\mathbb{R}^m$, and finally apply Lemma 1 in $\mathbb{R}^m.$ Again for details we refer to \cite{bg}. Finally we state the following lemma which follows from parabolic rescaling.
\begin{lemma}
Let $U_{\rho}$ be a cap of radius $\rho$ on S. Then
\[
\big\|\int_{U_{\rho}}f(\xi)e^{-ix\cdot \xi}d\sigma(\xi)\big\|_{L^p(B_R)} \leq {\rho}^{n-1-(n+1)/p}Q_{\rho R}^p.\]
\end{lemma}

Now we are in a position to describe the Bourgain-Guth argument. Let $f\in L^{\infty}(S), |f| \leq 1$ and let $x \in B_R$. Let 
\[ R^{\epsilon} \gg K_n \gg K_{n-1} \gg \ldots \gg K_1\]
 be constants independent of $f$ that will be  specified later. Decompose $S$ into caps $U^n_{\alpha}$ of size $1/K_n$, and let $\xi_{\alpha}^n \in U_{\alpha}^n$, thus
\begin{align*}
 \int_{S}f(\xi)e^{-ix\cdot \xi}\sigma(d\xi)&= \sum_{\alpha}e^{-ix \cdot \xi_{\alpha}^n}\int_{U_{\alpha}^n}f(\xi)e^{-ix\cdot (\xi-\xi_{\alpha}^n)}d\sigma(\xi)\\ &=: \sum_{\alpha}e^{-ix \cdot \xi_{\alpha}^n}T_{\alpha}^nf(x).
 \end{align*}
Take a function $\eta \in \mathcal{S}(\mathbb{R}^n)$  with 
$\widehat{\eta} (x)=1$ on $B(0,1)$ and $\widehat{\eta}(x)=0$ outside $B(0,2)$. Let $\eta_r(x)=\frac{1}{r^n}\eta(\frac{x}{r})$ and hence
\[
T_{\alpha}^nf(x)=T_{\alpha}^nf\ast \eta_{K_n}. \] 
For fixed $x$ using  the Bernstein inequality we can write
\begin{align*}
|T_{\alpha}^nf(x)| &\leq \int|T_{\alpha}^nf(x-y)\eta_{K_n}(y)|dy\\ 
&\lesssim \int|T_{\alpha}^nf(x-y)\eta_{K_n}(y)|^{\frac{1}{n}}|T_{\alpha}^nf(x-y)\eta_{K_n}(y)|^{\frac{n-1}{n}}dy \\
&\leq \big(\int |T_{\alpha}^nf(x-y)\eta_{K_n}(y)|^{\frac{1}{n}}dy\big) \|T_{\alpha}^nf(x-\cdot)\eta_{K_n}(\cdot)\|_{\infty}^{\frac{n-1}{n}}\\
&\lesssim \big(\int |T_{\alpha}^nf(x-y)\eta_{K_n}(y)|^{\frac{1}{n}}dy\big) \|T_{\alpha}^nf(x-\cdot)\eta_{K_n}(\cdot)\|_{1}^{\frac{n-1}{n}}K_n^{1-n}\\
\end{align*}
Thus if 
$|T_{\alpha}^nf(x)|$
is non-zero we can use division to write
\[\big(\int|T_{\alpha}^nf(x-y)\eta_{K_n}(y)|dy\big)^{1/n}\lesssim(\int |T_{\alpha}^nf(x-y)\eta_{K_n}(y)|^{\frac{1}{n}}dy)K_n^{1-n}.\]
 Hence
\[|T_{\alpha}^nf(x)|\lesssim\big(\int |T_{\alpha}^nf(x-y)\eta_{K_n}(y)|^{\frac{1}{n}}dy\big)^nK_n^{n-n^2}\]
This, of course, is trivial if $|T_{\alpha}^nf(x)|=0$, so it is independent of the value of $|T_{\alpha}^nf(x)|$. Taking the constant term inside the integral we obtain
\[\lesssim \big(\int |T_{\alpha}^nf(x-y)|^{1/n}\frac{1}{K^n_n}|\eta(\frac{y}{K_n})|^{1/n}dy \big)^{n}. \]
Define $\zeta(y)=\underset{|y-y'| \leq 1}{\max}|\eta(y')|^{1/n}$ and let $\zeta_r(y):=\frac{1}{r^n}\zeta(\frac{y}{r})$. Then 
\[|T_{\alpha}^nf(x)| \lesssim \big(\int|T_{\alpha}^nf(x-y)|^{1/n}\zeta_{K_n}(y)dy\big)^{n}=:c_{\alpha}^n(x)\] 
and
\[\big|\int_S f(\xi)e^{ix\cdot \xi}d\sigma(\xi) \big| \lesssim \sum_{\alpha}c_{\alpha}^n(x).\]
We have  $c_{\alpha}^n(x_1)\approx c_{\alpha}^n(x_2)$ whenever $|x_1-x_2|<K_n$.

 For a fixed $x$ we have two possibilities:

\textbf{1.1.} There exist $\alpha_1, \ldots, \alpha_n$ such that $|\xi_1' \wedge \ldots \wedge \xi_n'|> c(K_n)$ for $\xi_i \in U_{\alpha_i}^n$ and 
\[c_{\alpha_i}^n(x)> K_n^{-n}\underset{\alpha}{\max} \ c_{\alpha}^n(x).\]
 We can choose the same $\alpha_1, \ldots, \alpha_n$ for all $x$ in a ball of radius $K_n$ owing to the fact that $c_{\alpha}^n(x)$ are constant on balls of this size. 

\textbf{1.2.} The negation of this, namely there exist an $(n-1)$-dimensional  subspace $V_{n-1}$   such that 
\[c_{\alpha}^n(x) \leq K_n^{-n} \underset{\alpha}{\max} \  c_{\alpha}^n(x)\]
 if $ \text{dist} (U_{\alpha}^n,\widehat{V}_{n-1}) \gtrsim 1/K_n$ where $\widehat{V}$ stands for the image of $V \cap S^{n-1}$ under the Gauss map. Since  $c_{\alpha}^n(x)$ are essentially constant on  balls of radius $K_n$, on such balls we can take the    linear subspace $V_{n-1}$ to be the same for all $x \in  U_{\alpha}^n$.

If \textbf{1.1} holds then since the number of caps is comparable to $ K_n^{n-1}$
\[\big|\int_S f(\xi)e^{-ix\cdot \xi}\sigma(d\xi) \big| \lesssim K^{n-1}_n \underset{\alpha}{\max} \ c_{\alpha}^n(x)\lesssim K_n^{2n-1}(\prod_{i=1}^nc_{\alpha_i}^n(x))^{1/n}\] 
and thus letting $B_{1.1}$ denote $x\in B_R$ satisfying \textbf{1.1} 
\[
\int_{B_{1.1}}\big|\int_S f(\xi)e^{-ix\cdot \xi}\sigma(d\xi) \big|^{p}dx  \lesssim K^{(2n-1)p}_n \sum_{\alpha_1,\ldots,\alpha_n}\int_{B_R}(\prod_{i=1}^nc_{\alpha_i}^n(x))^{\frac{p}{n}}dx.\]
Then by definition of $c_{\alpha_i}^n(x)$ and H\"{o}lder inequality we have
\[\lesssim K_n^{(2n-1)p}\sum_{\alpha_1,\ldots,\alpha_n}\int_{B_R} \big( \prod_{i=1}^n\int|T_{\alpha}^nf(x-y_i)|^{p/n}\zeta_{K_n}(y_i)dy_i\big) dx\]
\[= K_n^{(2n-1)p}\sum_{\alpha_1,\ldots,\alpha_n}\int_{B_R} \big( \int \prod_{i=1}^n|T_{\alpha}^nf(x-y_i)|^{p/n}\zeta_{K_n}(y_i)dy_1\ldots dy_n\big) dx.\]
By Fubini's theorem 
\[\lesssim K_n^{(2n-1)p}\sum_{\alpha_1,\ldots,\alpha_n}\int \big( \int_{B_R} \prod_{i=1}^n|T_{\alpha}^nf(x-y_i)|^{p/n}dx\big)\prod_{i=1}^n\zeta_{K_n}(y_i)dy_1\ldots dy_n \]
Assuming $p\geq 2n/n-1$ the inner integral by \eqref{eq1.6'} satisfies 
\[\lesssim R^{\epsilon}\]
hence the main expression satisfies
\[\lesssim K_n^{10n^2}R^{\epsilon}\lesssim R^{2\epsilon}.\]
 The exponent $2n/n-1$ is the one prescribed by the restriction conjecture, thus in this $n$-linear case we get the best possible exponent.

Now assume 1.2 holds. Then since the number of caps is comparable to $K_n^{n-1}$ we have
\begin{align*}
\big| \int_S f(\xi)e^{-ix\cdot \xi}d\sigma(\xi) \big| & \lesssim \big| \int_{\{ \xi: \text{dist}(\xi,\hat{V}_{n-1})  \lesssim \frac{1}{K_n} \}} f(\xi)e^{-ix\cdot \xi}d\sigma(\xi) \big| \\ &+ \frac{1}{K_n} \underset{\alpha}{\max} \  c_{\alpha}^n(x)\\
&=\text{I} + \text{II}
\end{align*}
Thus
\begin{align*}
\int_{B_{1.2}}\big| \int_S f(\xi)e^{-ix\cdot \xi}d\sigma(\xi) \big|^p dx & \lesssim \int_{B_{1.2}} \big| \int_{\{ \xi: \text{dist}(\xi,\hat{V}_{n-1})  \lesssim \frac{1}{K_n} \}} f(\xi)e^{-ix\cdot \xi}d\sigma(\xi) \big|^p dx \\ &+ \frac{1}{K_n^p}\int_{B_R} \big(\underset{\alpha}{\max} \  c_{\alpha}^n(x)\big)^pdx
\end{align*}

We first evaluate the contribution coming from II. We have
\begin{align*}
 \frac{1}{K_n^p}\int_{B_R} (\underset{\alpha}{\max} \  c_{\alpha}^n(x))^p dx   & \leq \frac{1}{K_n^p} \sum_{\alpha} \int_{B_R} (c_{\alpha}^n(x))^p dx \\ &= \frac{1}{K_n^p}\sum_{\alpha}\| c_{\alpha}^n(x)\|^p_{L^p(B_R)}
\end{align*}
where the summation is over all caps of size $1/K_n$.   Using first H\"{o}lder  inequality then Lemma 3 we have,
\begin{align*}
\| c_{\alpha}^n(x)\|^p_{L^p(B_R)} & \lesssim \int_{B_R}\big(\int |T_{\alpha}^nf(x-y)|^p \zeta_{K_n}(y)dy\big)dx \\  & \lesssim \int \big(\int_{B_R} |T_{\alpha}^nf(x-y)|^p dx \big)\zeta_{K_n}(y) dy \\ & \lesssim K_n^{n+1-p(n-1)}Q^p_{R/K_n}               
\end{align*}
Hence contribution of II is bounded by
\[K_n^{n(2-p)}Q^p_{R/K_n}.\]
Since this is valid for all $f$ for an inductive argument aiming to bound $Q_R^p$ this is harmless when $p>2$.

To evaluate I we proceed as before, and decompose $S$ into caps $U^{n-1}_{\alpha}$ of radius $1/K_{n-1}$ and let $\xi^{n-1}_{\alpha} \in U^{n-1}_{\alpha}$. Then 

\begin{align*}
\int_{\{\xi: \text{dist}(\xi,\hat{V}_{n-1})  \lesssim \frac{1}{K_n}\}} f(\xi)e^{-ix\cdot      \xi}d\sigma(\xi) &= \sum_{\alpha} \int_{U_{\alpha}^{n-1} \cap \{\xi:\text{dist}(\xi,\hat{V}_{n-1})\lesssim \frac{1}{K_n}\}}f(\xi)e^{-ix\cdot      \xi}d\sigma(\xi)\\ &=:\sum_{\alpha}e^{-ix \cdot \xi_{\alpha}^{n-1}}\widetilde{T}_{\alpha}^{n-1}f(x)
\end{align*}
Let on the other hand
\[T^{n-1}_{\alpha}f(x)=\int_{U_{\alpha}^{n-1} }f(\xi)e^{-ix\cdot      \xi}d\sigma(\xi).\]
Note the small difference between  $\widetilde{T}$ and $T$, the first is defined on intersection of caps with a strip, while the second on full caps.
  Defining  $\eta_{K_{n-1}}$ and $\zeta_{K_{n-1}}$ just as $\eta_{K_n}$ and $\zeta_{K_n}$ were defined, we write
\[\widetilde{c}_{\alpha}^{n-1}(x):= \big(\int|\widetilde{T}_{\alpha}^{n-1}f(x-y)|^{1/n-1}\zeta_{K_{n-1}}(y)dy\big)^{n-1}.\]
\[c_{\alpha}^{n-1}(x):= \big(\int|T_{\alpha}^{n-1}f(x-y)|^{1/n-1}\zeta_{K_{n-1}}(y)dy\big)^{n-1}.\]
We shall need $c^{n-1}$ in the next step of our process. Via the same arguments as in in the definition of $c^n_{\alpha}$  we see that
\[|\widetilde{T}_{\alpha}^{n-1}f(x)|\lesssim \widetilde{c}_{\alpha}^{n-1}(x).\]
\[|T_{\alpha}^{n-1}f(x)|\lesssim c_{\alpha}^{n-1}(x).\]

We again have two cases for a fixed $x$:

\textbf{2.1.} There exist $\alpha_1,\ldots,\alpha_{n-1}$ such that $|\xi_1'\wedge \ldots \wedge \xi_{n-1}'|>c(K_{n-1})$ for $\xi_i \in U_{\alpha_i}^{n-1}$, and 
\[|\widetilde{c}_{\alpha_i}^{n-1}(x)| > K_{n-1}^{-(n-1)}\underset{\alpha}{\max} |\widetilde{c}_{\alpha}^{n-1}(x)|.\] 
Since $c_{\alpha}^{n-1}$ are essentially constant on balls of size $K_{n-1}$ we can choose $\alpha_{1},\ldots,\alpha_{n-1}$ the same for all $x$ in such balls.

\textbf{2.2.} The negation of this, namely there exist an (\textit{n-2})-dimensional  subspace $V_{n-2}$  which can be chosen a subspace of $V_{n-1}$ such that 
\[|\widetilde{c}_{\alpha}^{n-1}(x)| \leq K_{n-1}^{-(n-1)} \underset{\alpha}{\max}|\widetilde{c}_{\alpha}^{n-1}(x)|\] 
if $\text{dist}(U_{\alpha}^{n-1},\hat{V}_{n-2}) \ \gtrsim 1/K_{n-1} $. We can choose the linear subspace $V_{n-2}$ the same for all  $x$ in a  ball of size  $K_{n-2}.$

First assume \textbf{2.1} holds for a fixed $x$. Then
\[ \big| \int_{\{ \xi: \text{dist}(\xi,\hat{V}_{n-1})  \lesssim \frac{1}{K_n} \}} f(\xi)e^{-ix\cdot \xi}d\sigma(\xi) \big| \lesssim K_{n-1}^{2n-3} \big( \prod_{i=1}^{n-1}\widetilde{c}_{\alpha_i}^{n-1}(x) \big)^{1/{n-1}}.\]
If $p \geq 2(n-1)/(n-2)$ then we proceed to use the multilinear theory  of \cite{bct} as follows: 
\begin{equation*}\label{eq2.1}\int_{B_{2.1}}\big| \int_{\{ \xi: \text{dist}(\xi,\hat{V}_{n-1})  \lesssim \frac{1}{K_n} \}} f(\xi)e^{-ix\cdot \xi}d\sigma(\xi) \big|^p dx
\end{equation*}
 \[
 \lesssim C(K_{n-1})\int_{B_R} \big( \prod_{i=1}^{n-1}\widetilde{c}_{\alpha_i}^{n-1}(x) \big)^{p/{n-1}}dx  
 \]
where choices of $\widetilde{c}_{\alpha_i}^{n-1}$ of course depend on $x$. Clearly
\[\lesssim C(K_{n-1})\int_{B_R} \big( \prod_{i=1}^{n-1}\sum_{\alpha_i}c_{\alpha_i}^n(x) \big)^{p/{n-1}}dx \]
\[\lesssim C(K_{n})\sum_{\alpha_1,\ldots,\alpha_{n-1}}\int_{B_R} \big( \prod_{i=1}^{n-1}c_{\alpha_i}^{n}(x) \big)^{p/{n-1}}dx\]
and applying the same arguments as in  case \textbf{1.1} we obtain
\[\lesssim C(K_{n})R^{\epsilon} \lesssim R^{2\epsilon}\]
Thus we assume $p<2(n-1)/(n-2)$. In this case we have

\[ \oint_{B(a,K_n)}\big| \int_{\{ \xi: \text{dist}(\xi,\hat{V}_{n-1})  \lesssim \frac{1}{K_n} \}} f(\xi)e^{-ix\cdot \xi}d\sigma(\xi) \big|^p dx\]
\begin{equation}\label{new1} \leq K^{(2n-3)p}_{n-1} \sum_{\alpha_1,\ldots,\alpha_{n-1}}\oint_{B(a,K_n)}\big( \prod_{i=1}^{n-1}\widetilde{c}_{\alpha_i}^{n-1}(x)\big)^{p/n-1}dx\end{equation}
 where the subspace $V_{n-1}$ remains the same for all $x\in B(a,K_n)$. The choice of $\alpha_1,\ldots,\alpha_{n-1}$ remains the same only in balls of size $K_{n-1}$, but  since  the subspace is the same, caps $U_{\alpha}^{n-1}$  are always chosen from those intersecting the set 
 \[\{ \xi: \text{dist}(\xi,\hat{V}_{n-1})  \lesssim \frac{1}{K_n} \}.\]
 We will  exploit  multilinearity partially. Consider  an individual integral from \eqref{new1} above;  since $p<2(n-1)/(n-2)$ we have by H\"{o}lder inequality 
\begin{equation}\label{new3}\oint_{B(a,K_n)}\big( \prod_{i=1}^{n-1}\widetilde{c}_{\alpha_i}^{n-1}(x)\big)^{\frac{p}{n-1}}dx\end{equation}\[\lesssim \Big( \oint_{B(a,K_n)}\big( \prod_{i=1}^{n-1}\widetilde{c}_{\alpha_i}^{n-1}(x)\big)^{\frac{2}{n-2}}dx \Big)^{\frac{p(n-2)}{2(n-1)}} \]
which  by the definition of $\widetilde{c}^{n-1}_{\alpha}$ satisfies
\[\lesssim \Big( \oint_{B(a,K_n)}\big( \prod_{i=1}^{n-1}\int|\widetilde{T}_{\alpha}^{n-1}f(x-y_i)|^{\frac{1}{n-1}}\zeta_{K_{n-1}}(y_i)dy_i\big)^{\frac{2(n-1)}{n-2}}dx \Big)^{\frac{p(n-2)}{2(n-1)}}.\]
Using first H\"{o}lder inequality then Fubini's theorem we have
\[\lesssim \Big( \oint_{B(a,K_n)}\big( \prod_{i=1}^{n-1}\int|\widetilde{T}_{\alpha}^{n-1}f(x-y_i)|^{\frac{2}{n-2}}\zeta_{K_{n-1}}(y_i)dy_i\big)dx \Big)^{\frac{p(n-2)}{2(n-1)}}\]
\[\lesssim \Big(\int \big(\oint_{B(a,K_n)} \prod_{i=1}^{n-1}|\widetilde{T}_{\alpha}^{n-1}f(x-y_i)|^{\frac{2}{n-2}}dx\big)\prod_{i=1}^{n-1}\zeta_{K_{n-1}}(y_i)dy_1\ldots dy_{n-1}\ \Big)^{\frac{p(n-2)}{2(n-1)}}.\]
Now apply Lemma 2 to the inner integral to obtain 
\[\lesssim K_n^{\epsilon}\Big(\int \big(\oint_{B(a,K_n)} \prod_{i=1}^{n-1}(\sum_{\alpha_i}|T_{\alpha_i}^{n}f(x-y_i)|^2)^{\frac{1}{n-2}}dx\big)\prod_{i=1}^{n-1}\zeta_{K_{n-1}}(y_i)dy_1\ldots dy_{n-1}\ \Big)^{\frac{p(n-2)}{2(n-1)}}\] 
where the summation is  over all $\alpha_i $ such that $U_{\alpha_i}^n \subset U_{\alpha_i}^{n-1}$ and $U_{\alpha_i}^n\cap\hat{V}_{n-1} \neq \emptyset $. Since $p< 2(n-1)(n-2)$ by H\"{o}lder inequality and Fubini's theorem
 \begin{equation*}\label{new2} \lesssim K_n^{\epsilon+(n-2)(\frac{p}{2}-1)}\Big( \oint_{B(a,K_n)}  \big(\prod_{i=1}^{n-1}\int \big( \sum_{\alpha_i}|T_{\alpha_i}^{n}f(x-y_i)|^p\big) \zeta_{K_{n-1}}(y_i)dy_i\big)dx \Big)^{\frac{1}{n-1}}\end{equation*}
\[ \lesssim K_n^{\epsilon+(n-2)(\frac{p}{2}-1)}\Big( \oint_{B(a,K_n)}   \big( \sum_{\alpha}\int |T_{\alpha}^{n}f(x-y)|^p \zeta_{K_{n-1}}(y)dy\big)^{n-1}dx \Big)^{\frac{1}{n-1}}.\]
From the definition of $c^n_{\alpha}$ 
\[ \lesssim K_n^{\epsilon+(n-2)(\frac{p}{2}-1)}\Big( \oint_{B(a,K_n)}   \big( \sum_{\alpha}\int \big(c_{\alpha}^n(x-y)\big)^p \zeta_{K_{n-1}}(y)dy\big)^{n-1}dx \Big)^{\frac{1}{n-1}}\]
\[ \lesssim K_n^{\epsilon+(n-2)(\frac{p}{2}-1)}\Big( \oint_{B(a,K_n)}   \big( \sum_{\alpha} \big(c_{\alpha}^n(x)\big)^p \big)^{n-1}dx \Big)^{\frac{1}{n-1}}\]
where $\alpha$ in the summation is unrestricted. At this point we use the fact that the sum inside is constant:
\begin{equation}\label{new4} \lesssim K_n^{\epsilon+(n-2)(\frac{p}{2}-1)} \oint_{B(a,K_n)}    \sum_{\alpha} \big(c_{\alpha}^n(x)\big)^p dx, \end{equation}
Integrating both sides over $B_{2.1}$ to obtain
\begin{equation}\label{ex2}\int_{B_{2.1}}\big( \prod_{i=1}^{n-1}\widetilde{c}_{\alpha_i}^{n-1}(x)\big)^{\frac{p}{n-1}}dx \end{equation} \[\lesssim  K_n^{\epsilon+(n-2)(\frac{p}{2}-1)} \int_{B_R}   \big( \sum_{\alpha}\int |T_{\alpha}^{n}f(x-y)|^p \zeta_{K_n}(y)dy\big)dx \]

\[\lesssim  K_n^{\epsilon+(n-2)(\frac{p}{2}-1)}     \sum_{\alpha}\int \big(\int_{B_R} |T_{\alpha}^{n}f(x-y)|^p dx\big)\zeta_{K_n}(y)dy \]
Now we apply rescaling to obtain
\[\lesssim K_{n}^{(n-2)(\frac{1}{2}-\frac{1}{p})+\frac{n-1}{p}-(n-1)+\frac{n+1}{p}+\epsilon}Q^p_R\lesssim K_n^{\epsilon+\frac{n+2}{p}-\frac{n}{2}}Q^p_R\]
Thus we finally get
\[\int_{B_{2.1}}\big| \int_{\{ \xi: \text{dist}(\xi,\hat{V}_{n-1})  \lesssim \frac{1}{K_n} \}} f(\xi)e^{-ix\cdot \xi}d\sigma(\xi) \big|^p dx \lesssim C(K_{n-1})K_n^{\epsilon+\frac{n+2}{p}-\frac{n}{2}}Q^p_R\]
Suitably choosing $K_n$ with respect to $K_{n-1}$ shows that $p> 2(n+2)/n$ make this term acceptable. Thus we obtain the condition $p> \min(2n/(n-1),2(n+2)/n)$. Then we proceed similarly to handle \textbf{2.2}. We define $\widetilde{c}^{n-2}$ from intersections of caps with a $1/K_{n-1}$ neighborhood of our subspace  and $c^{n-2}$ from full caps. In our analysis $\widetilde{c}^{n-1}$ will then be replaced by $\widetilde{c}^{n-1}$, and $c^n$  by $c^{n-1}$. Continuing this process gives the condition
\begin{equation}\label{eq2.6}
p>2 \min(\frac{k}{k-1},\frac{2n-k+1}{2n-k-1}) \hspace{5mm} \text{for all} \hspace{5mm} 2\leq k \leq n
\end{equation}
which gives (\ref{eq1.4}). Here $p$ should be greater than  minimum for all values of $2\leq k \leq n$. Thus the condition is that $p$ should be greater than the maximum of these minima. When $n \equiv 0 \mod 3$  this maximum is attained only at the value $k=2n/3$ and it comes from the second term, thus has the value $2(4n+3)/(4n-3)$. The value of the first term, then, is $4n/(2n-3)$ which is strictly greater than  $\frac{8n+6}{4n-3}$. This is what allows us to improve the exponent when $n \equiv 0 \mod 3$. For the particular case of $n=6$ we have
\begin{equation}\label{eq2.5}
p>2 \min({\frac{k}{k-1},\frac{13-k}{11-k}}) \hspace{5mm} \text{for all}  \hspace{5mm}2\leq k \leq 6.
\end{equation}
For $k=6,5$ the minimum comes from the first term, for other $k$ from the second. The maximum of these minima comes from $k=4$ and is $18/7$. For  $k=4$ the first term gives $8/3>18/7$ and in our refined analysis we  shall exploit this.

\section{Refined Analysis}
From the analysis of the last section for a fixed $x\in B(0,R)$ we can write
\begin{equation}\label{eq3.1.1}
\begin{aligned}
|Tf(x)| & \leq C(K_6) \underset{i_1,\ldots,i_6}{\max}(c^6_{i_1}(x)\ldots c^6_{i_6}(x))^{1/6}\\ &+  \sum_{m=2}^5 C(K_m)\underset{V_{m}}{\max}\prod_{k=1}^m(\widetilde{c}^m_{i_k}(x))^{1/m}\\ &+  C(K_1)\underset{\alpha}{\max}\ c_{\alpha}^2(x) .
\end{aligned} 
\end{equation}
Here by $i_k$ we denote caps of radius $1/K_m$ in the $m$-linear term, $V_{m}$  denotes $m$ dimensional subspaces and $\hat{V}_{m}$ the image of $V_{m}\cap S^{5}$ under the Gauss map.  Our idea is, as among these terms $4$-linear gives the worst exponent, iterating this decomposition for smaller caps in that term we may obtain some improvement. We will also iterate the decomposition for linear, bilinear and trilinear terms. To execute this idea we replace the terms we want to further decompose as follows. Let $L_{m}$ denote the  caps $j$ of size $1/K_{m+1}$  such that $j \cap \hat{V}_{m} \neq \emptyset$. The  calculation we did to bound \eqref{new3} using Lemma 2 without the use of H\"{o}lder inequality  to raise the exponent   to $p$ in \eqref{new4} gives
\[\Big[\oint_{B(a,K_5)} \prod_{k=1}^4(\widetilde{c}^4_{i_k}(x))^{2/3}dx \Big]^{\frac{3}{8}}\lesssim K^{\epsilon}_5\big(\sum_{j\in L_4}(c_j^{5}(x))^2\big)^{1/2} \]
where on the right hand side $x\in B(a,K_5).$ 
Using this and the fact that $\widetilde{c}^4(x)$ are constant on balls of radius $K_4$ we  may write 
\begin{equation}\label{new5}
\prod_{k=1}^4(\widetilde{c}^4_{i_k}(x))^{1/4}=\phi_4 \cdot \big(\sum_{j\in L_4}(c_j^{5}(x))^2\big)^{1/2}
 \end{equation}
where $\phi_4 $ is constant on balls of radius $1$ and satisfies
\[\big( \oint_{ B(a,K_5)}(\phi_4)^{\frac{8}{3}}\big )^{\frac{3}{8}} \lesssim K^{\epsilon}_5.\]
Notice that if summation on the right hand side of \eqref{new5} is  zero then the left hand side is also zero, hence the function $\phi_4$ can  be constructed simply by dividing the left handside term by the summation on the right when summation is not zero, and by setting to zero when it is. 
For bilinear and trilinear terms we similarly have
\[\Big[\oint_{B(a,K_4)} \prod_{k=1}^3\widetilde{c}^3_{i_k}(x)dx \Big]^{\frac{1}{3}}\lesssim K^{\epsilon}_4\big(\sum_{j\in L_3}(c_j^{4}(x))^2\big)^{1/2} \]
\[\Big[\oint_{B(a,K_3)} \prod_{k=1}^2(\widetilde{c}^2_{i_k}(x))^{2}dx \Big]^{\frac{1}{4}}\lesssim K^{\epsilon}_3\big(\sum_{j\in L_2}(c_j^{3}(x))^2\big)^{1/2} .\]
So we can find $\phi_3$ and $\phi_2$ that are constant on balls of unit size that satisfy
 \[\prod_{k=1}^3(\widetilde{c}^3_{i_k}(x))^{1/3}=\phi_3 \cdot \big(\sum_{j\in L_3}(c_j^{4}(x))^2\big)^{1/2}\]
 \[  \prod_{k=1}^2(\widetilde{c}^2_{i_k}(x))^{1/2}=\phi_2 \cdot \big(\sum_{j\in L_2}(c_j^{3}(x))^2\big)^{1/2}  \]
\[\big( \oint_{ B(a,K_4)}(\phi_3)^3\big )^{\frac{1}{3}} \lesssim K^{\epsilon}_4, \ \ \  \ \big( \oint_{ B(a,K_3)}(\phi_2)^{4}\big )^{\frac{1}{4}} \lesssim K^{\epsilon}_3.\]
Thus we can write
\begin{align*}
|Tf(x)| & \leq C(K_6) \underset{i_1,\ldots,i_6}{\max}(c^6_{i_1}(x)\ldots c^6_{i_6}(x))^{1/6}\\ &+   C(K_5)\underset{V_{5}}{\max}\prod_{k=1}^5(\widetilde{c}^5_{i_k}(x))^{1/5}\\ &+  
\sum_{m=2}^4 K^{\epsilon}_{m+1} \phi_m \cdot \underset{V_{m}}{\max}\big(\sum_{j\in L_m}(c_j^{m+1})^2\big)^{1/2}\\ &+
C(K_1)\underset{\alpha}{\max}\ c_{\alpha}^2(x) .
\end{align*} 
 We shall iterate our decomposition for these caps $j$.  We now describe this in a general fashion. Let $\tau$ be a cap of radius $\delta$. By first scaling to the unit scale, then applying the decomposition, then scaling back we get the following:
\begin{align*}
|Tf(x)| & \leq C(K_6) \underset{i_1,\ldots,i_6}{\max}(c^6_{\tau_1}(x)\ldots c^6_{\tau_6}(x))^{1/6}\\ &+   C(K_5)\underset{V_{5}}{\max}\prod_{k=1}^5(\widetilde{c}^5_{\tau_k}(x))^{1/5}\\&+
\ \sum_{m=2}^4K_{m+1}^{\epsilon} \phi_{\tau_m} \cdot\underset{V_{m}}{\max} \big(\sum_{j\in L_m}(c_{\eta}^{m+1}(x))^2\big)^{1/2}\\ &
 +  C(K_1)\underset{\alpha}{\max}\ c_{\tau_{\alpha}}^2(x) .
\end{align*}
Here, similar to what we have above by $\tau_k$ we denote caps that are of radius $1/K_m$ in $m$-linear term, and the notation $\eta$ denotes caps of radius $1/K_{m+1}$ in $m$-linear term. Furthermore we have $\phi_{\tau_m}$ constant on boxes $\tau'$ that  are dual to cap $\tau$ and for boxes $K_{m+1}\tau'$ 
\[\oint_B\phi_{\tau_m}^{8/3}\lesssim K_{m+1}^{\epsilon}.\]
These two are simple consequences of rescaling. 

We iterate this decomposition  except for 6-linear and 5-linear terms in  our main decomposition. We clarify  several points that arises from application of this process. First of all from second step onwards we actually apply the decomposition not to expressions of type $Tf_{\tau}$ but $c_{\tau}$. This is a simple issue to deal with, and right hand side remains  the same. To see this notice that all terms on the right hand side are already constant on balls exceeding the size of averaging we need to pass from $T_{\tau}$ to $c_{\tau}$.   Secondly,  as we iterate, we will need to multiply functions $\phi_{\tau_m}$ arising in each step of iteration. To investigate what happens in this case let  $\phi_{\tau_k}$ and 
$\phi_{\eta_l}$ be such functions arising in consecutive steps.  Thus $\phi_{\tau_k}$ is constant on boxes $\tau'$  that are dual to the cap $\tau$, and  $\phi_{\eta}$
constant on boxes $\eta'$ that are dual to the cap $\eta$. As $\eta$ arise when we decompose $\tau$, if we let $\tau$ be a $\delta$ cap, $\eta$ is a $\delta/K_{l+1}$ cap.
Now let $B$ be a $K_{l+1}\eta'$ box. One can, of course, decompose this box into $\eta'$ boxes $B_{\alpha}$; thus since  $\phi_{\eta_l}^{8/3}$ is comparable to a constant on a $B_{\alpha}$ box
\[\int_B \phi_{\tau_k}^{8/3}\phi_{\eta_l}^{8/3} \lesssim \sum_{\alpha}\phi_{\eta_l}^{8/3}\Big|_{B_{\alpha}}\int_{B_{\alpha}} \phi_{\tau_k}^{8/3}\]
where  we  use the expression 
\[\phi_{\eta_l}^{8/3}\Big|_{B_{\alpha}}\]
 to denote evaluation of $\phi_{\eta_l}^{8/3}$ on an arbitrary point of $B_{\alpha}$.
The direction of a $\tau'$ box differs from that of a $\eta'$ only by an angle of $\delta$, thus decomposing it we should have average over the larger box at most maximum of averages over smaller boxes. Hence
 \[ \oint_{B_{\alpha}} \phi_{\tau_k}^{8/3} \ \lesssim K_{k+1}^{\epsilon}.\]
So we have
\begin{align*}
\sum_{\alpha}\phi_{\eta_l}^{8/3}\Big|_{B_{\alpha}}\int_{B_{\alpha}} \phi_{\tau_k}^{8/3}& \lesssim \sum_{\alpha}\int_{B_{\alpha}}\phi_{\eta_l}^{8/3}\oint_{B_{\alpha}} \phi_{\tau_k}^{8/3} \\ &\lesssim K^{\epsilon}_{k+1} \sum_{\alpha} \int_{B_{\alpha}} \phi_{\eta_l}^{8/3}\\
 & \lesssim K_{k+1}^{\epsilon} \int_B  \phi_{\eta_l}^{8/3}\\ & \lesssim K_{k+1}^{\epsilon}K_{l+1}^{\epsilon}|B|. 
\end{align*}
This also shows that emerging cross terms will not lead to any problem in the iteration process as losses $K_{m+1}^{\epsilon}$ are proportionate to size of caps $\tau_m$. That is if we divide a cap into larger caps, the number of steps we iterate our decomposition  increases, but since the losses incurred at each step are smaller this does not lead to any problem.
After this investigation we are ready to state the final situation after iterating the decomposition:
\begin{align*}
|Tf| &\leq R^{\epsilon} \underset{R^{-1/2}<\delta <1}{\max} \underset{E_{\delta}}{\max} \Big[ \sum_{\tau \in E_{\delta}}\big(\phi_{\tau}\underset{\tau_1,\ldots,\tau_6}{\max}(c^6_{\tau_1}\ldots c^6_{\tau_6})^{1/6}  \big)^2 \Big]^{1/2} \\ 
&+ R^{\epsilon}   \underset{R^{-1/2}<\delta <1}{\max} \underset{E_{\delta}}{\max} \Big[\sum_{\tau \in E_{\delta}} \big(\phi_{\tau} (\underset{V_{4}}{\max} \prod_{k=1}^{5}  \widetilde{c}^{5}_{\tau_k})^{1/5}\big)^2  \Big]^{1/2} \\
&+ R^{\epsilon} \underset{E_{R^{-1/2}}}{\max} \Big[ \sum_{\tau \in E}(\phi_{\tau}c_{\tau})^2 \Big]^{1/2}
\end{align*}
where 

\textbf{1.} In all terms $E_\delta$ is a collection of $\delta$ caps, with cardinality $\delta^{-3}$.

\textbf{2.} In the $m$-linear term $\tau_j \subset \tau$ are caps of size $\delta / K_m$ satisfying $\tau_j \cap \hat{V}_{m-1} \neq \emptyset$ and the linear independence condition.

\textbf{3.}  For  $B$  a $\tau'$ box we have
\[\oint_B\phi_{\tau}^{8/3} < R^{\epsilon}.\]
We shall estimate each of the following terms above. For each of the $m$-linear terms with $m>1$ we proceed as follows: we first estimate the term in the $L^p$ space with  the exponent given by (\ref{eq2.5}) for $k=m$.  Then we estimate at $L^{18/7}$  which is the maximum of  exponents given by  (\ref{eq2.5}) exactly in the same fashion. Using this, we do a refined estimate using pigeonholing at the exponent $18/7$ in two different ways to obtain a small gain. Interpolation with the estimate at the exponent given by (\ref{eq2.5}) will determine the small amount of improvement to the exponent $18/7$. For the linear term the process is similar but simpler. We start with estimating the six-linear term.

\subsection{Estimates on the six-linear term}
We consider the term
\begin{equation}\label{eq3.1.1}
\underset{E_{\delta}}{\max} \Big[ \sum_{\tau \in E_{\delta}}\big(\phi_{\tau}\underset{\tau_1,\ldots,\tau_6}{\max}(c^6_{\tau_1}\ldots c^6_{\tau_6})^{1/6}  \big)^2 \Big]^{1/2}.
\end{equation}
The inner maximum has no importance and we can fix $c^6_{\tau_i}$. We wish to exploit the fact that Fourier transform of a function supported on a small cap is roughly constant on tubes dual to this cap. To make this precise we shall need some further notation. Recall that we took a function  $\eta \in \mathcal{S}(\mathbb{R}^n)$  with 
$\widehat{\eta} (x)=1$ on $B(0,1)$ and $\widehat{\eta}(x)=0$ outside $B(0,2)$. Rescale this function to obtain $\upsilon_{\tau_i}$ adapted  to  tubes dual to the cap $\tau_i$. Similarly obtain  $\beta_{\tau_i}$ by rescaling $\zeta$.
With $\upsilon_{\tau_i}, \beta_{\tau_i}$ we define $b_{\tau_i}^6$ as $c_{\tau_i}^6$ is defined in page 5. Then 
\[|Tf_{\tau_i}| \lesssim b_{\tau_i}^6. \]
 So convolving both sides with $\zeta_{\tau_i}$ and using the fact that $b_{\tau_i}^6$ are actually constant on balls of size that this averaging takes place we obtain
 \[c_{\tau_i}^6 \lesssim b_{\tau_i}^6. \]
Thus we can estimate (\ref{eq3.1.1}) by 
\[\underset{E_{\delta}}{\max} \Big[ \sum_{\tau \in E_{\delta}}\big(\phi_{\tau}(b^6_{\tau_1}\ldots b^6_{\tau_6})^{1/6}  \big)^2 \Big]^{1/2}.\]
Assume $|f| \leq 1$. We have
\[\int_{B_R}(b^6_{\tau_1}\ldots b^6_{\tau_6})^{2/5} =\int_{B_R}\big(\prod_{i=1}^6\int |Tf_{\tau_i}(x-y)|^{1/6}_{\delta/K_6}\beta_{\delta/K_6}(y)dy \big)^{12/5}dx   \]
by H\"{o}lder inequality 
\[\lesssim \int_{B_R}\big(\prod_{i=1}^6\int |Tf_{\tau_i}(x-y)|^{2/5}\beta_{\delta/K_6}(y)dy \big)dx \]
\[= \int_{B_R}\big(\int \prod_{i=1}^6 |Tf_{\tau_i}(x-y_i)|^{2/5}\beta_{\delta/K_6}(y_i)dy_1\ldots dy_6 \big)dx\]
which by Fubini's theorem becomes
\[=\int \big(\int_{B_R} \prod_{i=1}^6 |Tf_{y_i,\tau_i}(x)|^{2/5}dx\big )\prod_{i=1}^6\beta_{\delta/K_6}(y_i)dy_1\ldots dy_6 \]
Of course $f_{y_i}$ are modulations of $f$. Now rescaling inside integral  to obtain functions $|g_{y_i}| \leq 1$ and caps  $U_1, \cdots , U_6$ of size $1/K_6$ satisfying the linear independence condition we have
\[\lesssim \delta^5\int \big(\int_{B_R} \prod_{i=1}^6 |Tg_{y_i,U_i}(x)|^{2/5}dx\big )\prod_{i=1}^6\beta_{\delta/K_6}(y_i)dy_1\ldots dy_6 \]
 Thus multilinear theory of \cite{bct} applies to yield
\[ \lesssim \delta^5 R^{\epsilon}.\]

With this at hand we proceed to estimate at the exponent $12/5$ given by (\ref{eq2.5}). Using H\"{o}lder we have
\begin{align*}  \Big[ \sum_{\tau \in E_{\delta}}\big(\phi_{\tau}(b^6_{\tau_1}\ldots b^6_{\tau_6})^{1/6}  \big)^2 \Big]^{1/2} 
&\leq |E_{\delta}|^{\frac{1}{12}} \Big[ \sum_{\tau \in E_{\delta}}\big(\phi_{\tau}(b^6_{\tau_1}\ldots b^6_{\tau_6})^{1/6}  \big)^{\frac{12}{5}} \Big]^{\frac{5}{12}} \\ &\lesssim  \delta^{-\frac{1}{4}}\Big[ \sum_{\tau \in E_{\delta}}\big(\phi_{\tau}(b^6_{\tau_1}\ldots b^6_{\tau_6})^{1/6}  \big)^{\frac{12}{5}} \Big]^{\frac{5}{12}}   
\end{align*} 
Now $\tau$ ranges over a full partition into $\delta$- caps of $S$, and does not depend on particular choice of $E_{\delta}$. Now let $B$ stand for $\tau'$ boxes. Since $b_{\tau_i}^6$ are constant on $\tau_i'$ boxes, we have 
\begin{align*}
\int_{B_R} \big(\phi_{\tau}(b^6_{\tau_1}\ldots b^6_{\tau_6})^{1/6}  \big)^{\frac{12}{5}}& 
\lesssim \sum_{B} (b^6_{\tau_1}\ldots b^6_{\tau_6})^{\frac{2}{5}} \Big|_B \big( \int_B \phi^{12/5}_{\tau}) \\ &
 \lesssim  \sum_{B} (\int _B(b^6_{\tau_1}\ldots b^6_{\tau_6})^{\frac{2}{5}}\big)  (\oint_B \phi^{12/5}_{\tau}) \\
& \lesssim R^{\epsilon} \int _{B_R}(b^6_{\tau_1}\ldots b^6_{\tau_6})^{\frac{2}{5}} \\
 & \lesssim R^{\epsilon}\delta^{5}.
\end{align*}
Using this we finally obtain
\begin{equation*}\label{eq3.1.2}
\|(\ref{eq3.1.1})\|_{L^{12/5}(B_R)} \lesssim R^{\epsilon}\delta^{-1/4}
.\end{equation*}
On the other hand applying the same process and using the fact  $b^6_{\tau_i}\lesssim \delta^5$ to reduce the exponent yields
\begin{equation*}\label{eq3.1.3}
\|(\ref{eq3.1.1})\|_{L^{18/7}(B_R)} \lesssim R^{\epsilon}.
\end{equation*}

Now we begin finer estimates. Let $0<\lambda<1$  and define
\[g_{\tau, \lambda} =g_{\tau}1_{\{g_{\tau} \sim \lambda \delta^5\}} \hspace{5mm} \text{where} \hspace{5mm} g_{\tau}=(b^6_{\tau_1}\ldots b^6_{\tau_6})^{1/6} .\]
Then 
\[\int_{B_R}g_{\tau,\lambda}^{18/7} < (\lambda \delta^5)^{18/7-12/5}\int_{B_R}g_{\tau,\lambda}^{12/5} \lesssim R^{\epsilon} \lambda^{6/35}\delta^{41/7}.  \]
Using this
\begin{equation*}\label{eq3.1.3}
\Big[ \int_{B_R} \underset{E_{\delta}}{\max} \big( \sum_{\tau \in E_{\delta} }(\phi_{\tau}g_{\tau,\lambda})^2  \big)^{9/7} \Big]^{7/18} \lesssim R^{\epsilon} \lambda^{1/15}.
\end{equation*}
We do one more pigeonholing. Let $1\leq \mu < \infty$ and decompose 
\[\phi_{\tau}=\sum_{\mu \hspace{1mm} \text{dyadic}}\phi_{\tau,\mu}\]
where
\[\phi_{\tau,\mu}=\phi_{\tau}1_{\{\phi_{\tau} \sim \mu\}}, \hspace{5mm} \phi_{\tau,1}=\phi_{\tau}1_{\{\phi_{\tau} \leq 1 \}} .\]
Then we have
\[ \oint_B \phi_{\tau,\mu}^{18/7} \leq \mu^{-2/21}\oint_B \phi_{\tau,\mu}^{8/3} \lesssim R^{\epsilon}\mu^{-2/21}. \]
\begin{equation*}\label{eq3.1.4}
\Big[ \int_{B_R} \underset{E_{\delta}}{\max} \big( \sum_{\tau \in E_{\delta} }(\phi_{\tau,\mu}g_{\tau,\lambda})^2  \big)^{9/7} \Big]^{7/18} \lesssim R^{\epsilon} \lambda^{1/15}\mu^{-1/27}.
\end{equation*}

We  now estimate the left hand side of the inequality above in a different way. Clearly we have
\begin{equation*}\label{eq3.5}
\underset{E_{\delta}}{\max} \big( \sum_{\tau \in E_{\delta} }(\phi_{\tau,\mu}g_{\tau,\lambda})^2  \big)^{1/2} \leq \mu \big( \sum_{\tau  }g_{\tau,\lambda}^2  \big)^{1/2}
.\end{equation*}
Now $\tau$ ranges over a full partition into $\delta$ caps of the surface $S$. We shall write the right hand side as convolutions of measures with tubes, and apply Kakeya maximal function estimates. 
Since separation between directions of caps $\tau$ and $\tau_i$ are small we have
\begin{equation}\label{new6}
(b_{\tau_i}^6)^{1/6} \lesssim (b_{\tau_i}^6)^{1/6} \ast 
(\delta^71_{\tau'}).
\end{equation}
Hence  
\begin{align*}
g_{\tau} & \lesssim \int \big( \prod_{i=1}^{6}(b_{\tau_i}^6)^{1/6}\ast \delta^71_{\tau'} \big)(z) (\delta^7 1_{\tau'})(x-z) dz \\
& \lesssim \int \omega (z) (\delta^{7} 1_{\tau'})(x-z)dz.
\end{align*}
We by the definition of $g_{\tau,\lambda}$ and \eqref{new6} have
\[g^2_{\tau,\lambda}\lesssim \omega^21_{\{\omega \gtrsim \lambda\delta^5\}}\]
But since the function $\omega$ is constant on tubes dual to $\tau$ we can write
\[g_{\tau,\lambda}^2 \lesssim \delta^{7}\int (\omega^21_{\{ \omega \gtrsim  \lambda    \delta^5 \} })(z) 1_{\tau'}(x-z)dz. \]
We wish to replace the expression 
 \[(\omega^21_{\{ \omega \gtrsim  \lambda    \delta^5 \} })(z)dz\] 
with a constant multiple of a probability measure, from which we will pass to the Kakeya maximal function. To this end we estimate the total mass of this measure. By Chebyshev's inequality
\[
\int_{B_R}\omega^21_{\{ \omega \gtrsim \lambda \delta^5 \} }(x)dx  \lesssim (\frac{1}{\lambda \delta^5})^{2/5}\int_{B_R} \omega^{12/5}(x)dx \]
which is
 \[ \lesssim (\frac{1}{\lambda \delta^5})^{2/5} \int \big(\int_{B_R} (\prod_{i=1}^6b_{\tau_i}^6(x-z_i))^{2/5}dx \big) \big(\prod_{i=1}^6 (\delta^71_{\tau'})(z_i) \big)dz_1\cdots dz_6
\]
\[\lesssim R^{\epsilon} \lambda^{-2/5} \delta^{3}. \]

This puts us in a position to bring into  play the Kakeya maximal function estimate of Wolff, for we have
\[g_{\tau,\lambda}^2 \lesssim R^{\epsilon}\delta^{10}\lambda^{-2/5} \int 1_{\tau'}(x-y)d\mu_{\tau}\]
where $d\mu_{\tau}$ is a probability distribution. From this and convexity we have
\begin{align*} 
\|\underset{E_{\delta}}{\max} \big( \sum_{\tau \in E_{\delta} }(\phi_{\tau,\mu}g_{\tau,\lambda})^2  \big)^{1/2}\|_{L^{18/7}(B_R)} & \leq R^{\epsilon}\mu \delta^5 \lambda^{-1/5} \| \big( \sum_{\tau}1_{\tau'}(x-y_{\tau}) \big)^{1/2} \|_{L^{18/7(B_R)}}\\ &= R^{\epsilon}\mu \delta^5 \lambda^{-1/5} \| \big( \sum_{\tau}1_{\tau'}(x-y_{\tau}) \big) \|_{L^{9/7(B_R)}}^{1/2}
 \end{align*}
 where $y_{\tau}$ is a choice of points in $\textbf{R}^6$. Now we are ready to apply Wolff's Kakeya estimate from \cite{w}. For $\delta$-separated $\delta$-tubes $T$ this estimate gives
 \[ \|\sum_{T}\chi_{T}\|_{L^{4/3}} \lesssim \delta^{-1/2-} .\]
 Interpolating this with the trivial estimate at $L^1$ yields 
 \[ \|\sum_{T}\chi_{T}\|_{L^{9/7}} \lesssim \delta^{-4/9-}. \]
 We rescale the  estimate since sizes of our tubes are different. This gives for our tubes
 \[ \|  \sum_{\tau}1_{\tau'}(x-y_{\tau})  \|_{L^{9/7}(B_R)}^{1/2} \lesssim  \delta^{-44/9}.\]
 Thus finally
 \[\|\underset{E_{\delta}}{\max} \big( \sum_{\tau \in E_{\delta} }(\phi_{\tau,\mu}g_{\tau,\lambda})^2  \big)^{1/2}\|_{L^{18/7}(B_R)} \lesssim R^{\epsilon}\mu \delta^{1/9} \lambda^{-1/5}\lesssim R^{\epsilon}\mu \delta^{1/9} \lambda^{-9/5}.\]

 Now let's see our estimates together:
\[\| \underset{E_{\delta}}{\max} \big( \sum_{\tau \in E_{\delta} }(\phi_{\tau}g_{\tau})^2  \big)^{1/2}\|_{L^{12/5}(B_R)} \lesssim R^{\epsilon}\delta^{-1/4}.\]

\begin{align*}\| \underset{E_{\delta}}{\max} \big( \sum_{\tau \in E_{\delta} }(\phi_{\tau,\mu}g_{\tau,\lambda})^2  \big)^{1/2}\|_{L^{18/7}(B_R)} & \lesssim R^{\epsilon}\min(\mu  \lambda^{-9/5}\delta^{1/9}, \lambda^{1/15}\mu^{-1/27})\\ &\lesssim R^{\epsilon}\delta^{\frac{1}{9\cdot 28}}.\end{align*}
Interpolating these gives the small improvement
\[\frac{2}{735}.\]

\subsection{Estimates on the $5$-linear term}

We consider the term 
\begin{equation}\label{eq3.2.1}
\underset{E_{\delta}}{\max} \Big[\sum_{\tau \in E_{\delta}} \big(\phi_{\tau} (\underset{V_{4}}{\max} \prod_{k=1}^{5}\widetilde{c}^{5}_{\tau_k})^{1/5}\big)^2  \Big]^{1/2}.
\end{equation}
 This, as above, can be estimated by
\[\underset{E_{\delta}}{\max} \Big[\sum_{\tau \in E_{\delta}} \big(\phi_{\tau} (\underset{V_{4}}{\max} \prod_{k=1}^{5}\widetilde{b}^{5}_{\tau_k})^{1/5}\big)^2  \Big]^{1/2}\]
The  maximum taken over all $V_{4}$ does not make  any difference to our estimates  since our estimates will remain the same for all $V_{4}$. Thus  fix $V_{4}$. Assume $|f|\leq1$. Using rescaling as in six-linear case one obtains
\begin{equation}\label{new7}
\int_{B_R} (\prod_{k=1}^{5}\widetilde{b}^{5}_{\tau_k})^{1/2} \leq \delta^{11/2}R^{\epsilon}.  \end{equation}
With this we proceed as in six-linear case. By H\"{o}lder 
\begin{align*} \Big[ \sum_{\tau \in E_{\delta}}\big(\phi_{\tau}( \prod_{k=1}^{5}\widetilde{b}^{5}_{\tau_k})^{\frac{1}{5}}\big)^2 \Big]^{\frac{1}{2}}  &\leq |E_{\delta}|^{\frac{1}{10}}\Big[ \sum_{\tau  }\big(\phi_{\tau}( \prod_{k=1}^{5}\widetilde{b}^{5}_{\tau_k})^{\frac{1}{5}}\big)^{\frac{5}{2}} \Big]^{\frac{2}{5}} \\ &\leq  \delta^{-\frac{3}{10}}\Big[ \sum_{\tau  }\big(\phi_{\tau}( \prod_{k=1}^{5}\widetilde{b}^{5}_{\tau_k})^{\frac{1}{5}}\big)^{\frac{5}{2}} \Big]^{\frac{2}{5}}.
\end{align*} 
Now $\tau$ ranges over a full partition into $\delta$- caps of $S$, and does not depend on particular choice of $E_{\delta}$. Now let $B$ stand for $\tau'$ boxes. Since $\widetilde{b}^{5}_{\tau_k}$ are constant on $\tau_i'$ boxes, we have 
\begin{align*}
\int_{B_R} \big(\phi_{\tau}(\prod_{k=1}^{5}\widetilde{b}_{\tau_k}^{5}\big)^{\frac{1}{2}}& \lesssim \sum_{B} ( \prod_{k=1}^{5}\widetilde{b}_{\tau_k}^{5})^{\frac{1}{2}}\Big|_B ( \int_B \phi^{5/2}_{\tau}) \\ & \lesssim  \sum_{B} (\int _B( \prod_{k=1}^{5}\widetilde{b}_{\tau_k}^{5})^{\frac{1}{2}}\big)  (\oint_B \phi^{5/2}_{\tau}) \\
& \lesssim R^{\epsilon} \int _{B_R}( \prod_{k=1}^{5}\widetilde{b}_{\tau_k}^{5})^{\frac{1}{2}} \\
 & \lesssim R^{\epsilon}\delta^{11/2}.
\end{align*}
Using this we finally obtain
\begin{equation*}
\|(\ref{eq3.2.1})\|_{L^{5/2}(B_R)} \lesssim R^{\epsilon}\delta^{-1/10}
.\end{equation*}
On the other hand the same  process yields
\begin{equation*}
\|(\ref{eq3.2.1})\|_{L^{18/7}(B_R)} \lesssim R^{\epsilon}
.\end{equation*}

Now we proceed to finer estimates.   Let this time
\[g_{\tau, \lambda} =g_{\tau}1_{\{g_{\tau} \sim \lambda \delta^5\}} \hspace{5mm} \text{where} \hspace{5mm} g_{\tau}=(\prod_{k=1}^{5}\widetilde{b}^5_{\tau_k})^{1/5}.\]
Then by \eqref{new7}
\[\int_{B_R}g_{\tau,\lambda}^{18/7} < (\lambda \delta^5)^{18/7-5/2}\int_{B_R}g_{\tau,\lambda}^{5/2} \lesssim R^{\epsilon} \lambda^{1/14}\delta^{41/7}.  \]
Using this
\begin{equation*}\label{eq3.1.3}
\Big[ \int_{B_R} \underset{E_{\delta}}{\max} \big( \sum_{\tau \in E_{\delta} }(\phi_{\tau}g_{\tau,\lambda})^2  \big)^{9/7} \Big]^{7/18} \lesssim R^{\epsilon} \lambda^{1/36}.
\end{equation*}
Decompose $\phi_{\tau}$ exactly as before. Then we have
\begin{equation*}\label{eq3.1.4}
\Big[ \int_{B_R} \underset{E_{\delta}}{\max} \big( \sum_{\tau \in E_{\delta} }(\phi_{\tau,\mu}g_{\tau,\lambda})^2  \big)^{9/7} \Big]^{7/18} \lesssim R^{\epsilon} \lambda^{1/36}\mu^{-1/27}.
\end{equation*}

Now we estimate 
\[\underset{E_{\delta}}{\max} \big( \sum_{\tau \in E_{\delta} }(\phi_{\tau,\mu}g_{\tau,\lambda})^2  \big)^{1/2}\]
using the Kakeya maximal function as in the six-linear case. Clearly
\begin{equation*}\label{eq3.5}
\underset{E_{\delta}}{\max} \big( \sum_{\tau \in E_{\delta} }(\phi_{\tau,\mu}g_{\tau,\lambda})^2  \big)^{1/2} \leq \mu \big( \sum_{\tau  }g_{\tau,\lambda}^2  \big)^{1/2}
\end{equation*}
where $\tau$ ranges over a full partition into $\delta$ caps of the surface $S$. \[(\widetilde{b}_{\tau_k}^5)^{1/5} \lesssim (\widetilde{b}_{\tau_k}^5)^{1/5} \ast 
(\delta^71_{\tau'}).\]
Hence  
\begin{align*}
g_{\tau} & \lesssim \int \big( \prod_{k=1}^{5}(\widetilde{b}_{\tau_i}^5)^{1/5}\ast \delta^71_{\tau'} \big)(z) (\delta^7 1_{\tau'})(x-z) dz \\
& \lesssim \int \omega (z) (\delta^{7} 1_{\tau'})(x-z)dz.
\end{align*}
So
\[g_{\tau,\lambda}^2 \lesssim \delta^{7}\int (\omega^21_{\{ \omega \gtrsim  \lambda    \delta^5 \} })(z) 1_{\tau'}(x-z)dz. \]
We  have
\[
\int_{B_R}\omega^21_{\{ \omega \gtrsim \lambda \delta^5 \} }(x)dx  \lesssim (\frac{1}{\lambda \delta^5})^{1/2}\int_{B_R} \omega^{5/2}(x)dx \]
which by \eqref{new7}
 \[ \lesssim (\frac{1}{\lambda \delta^5})^{1/2} \int_{B_R} \big(\int (\prod_{k=1}^5\widetilde{b}_{\tau_i}^5(x-z_i))^{\frac{1}{2}}dx \big) \big(\prod_{k=1}^5 (\delta^71_{\tau'})(z_i) \big)dz_1\cdots dz_5
\]
\[\lesssim R^{\epsilon} \lambda^{-1/2} \delta^{3}. \]
Hence,
\[g_{\tau,\lambda}^2 \lesssim R^{\epsilon}\delta^{10}\lambda^{-1/2} \int 1_{\tau'}(x-y)d\mu_{\tau}(dy).\]
Thus proceeding exactly as in the six-linear case we shall have 
 \[\|\underset{E_{\delta}}{\max} \big( \sum_{\tau \in E_{\delta} }(\phi_{\tau,\mu}g_{\tau,\lambda})^2  \big)^{1/2}\|_{L^{18/7}(B_R)} \lesssim R^{\epsilon}\mu \delta^{1/9} \lambda^{-1/4}\lesssim R^{\epsilon}\mu \delta^{1/9} \lambda^{-3/4}.\]
 Summarizing our estimates
\[\| \underset{E_{\delta}}{\max} \big( \sum_{\tau \in E_{\delta} }(\phi_{\tau}g_{\tau})^2  \big)^{1/2}\|_{L^{5/2}(B_R)} \lesssim R^{\epsilon}\delta^{-1/10}.\]

\begin{align*}\| \underset{E_{\delta}}{\max} \big( \sum_{\tau \in E_{\delta} }(\phi_{\tau,\mu}g_{\tau,\lambda})^2  \big)^{1/2}\|_{L^{18/7}(B_R)} & \lesssim R^{\epsilon}\min(\mu  \lambda^{-3/4}\delta^{1/9}, \lambda^{1/36}\mu^{-1/27})\\ &\lesssim R^{\epsilon}\delta^{\frac{1}{9\cdot 28}}.\end{align*}
Then interpolation yields that improvement is 
\[\frac{5}{1764}.\]

\subsection{Estimates on the linear term}
 For this term caps are of size  $R^{-1/2}$ thus we can use more direct methods without suffering any significant loss. One such method is using the fact that for a function  $f$ supported on $S$ we have for every $x_n\in \mathbb{R}^n$
 \[\|\widehat{fd\sigma}\|_{L^2(\mathbb{R}^{n-1}\times \{x_n\})} \approx \|f\|_{L^2(S)}\]
 This is referred to as conservation of mass in PDE literature. We shortly give the calculation that leads to this result. Let $x=(\overline{x},x_n).$ Let $S$ be parametrized by $\xi_n=\phi(\overline{\xi})$. Then 
 \begin{align*}\|\widehat{fd\sigma}\|_{L^2(\mathbb{R}^{n-1}\times \{x_n\})}^2&=\int|\int f(\overline{\xi})e^{-i(\overline{x}\cdot\overline{\xi} +x_n\phi(\overline{\xi})}d\overline{\xi}|^2d\overline{x} \\
&=\int|\int f(\overline{\xi})e^{x_n\phi(\overline{\xi})}e^{-i\overline{x}\cdot\overline{\xi} }d\overline{\xi}|^2d\overline{x}.
\end{align*} 
Now the inner integral is the Fourier transform of $f(\overline{\xi})e^{x_n\phi(\overline{\xi})}$, thus applying the Plancherel theorem we have
\[=\int|f(\overline{\xi})e^{x_n\phi(\overline{\xi})}|^2d\overline{\xi}=\int|f(\overline{\xi})|^2d\overline{\xi} \approx \|f\|_{L^2(S)}^2.\]
  
  To continue estimation  of the linear term  we pass from $c_{\tau}$  to $b_{\tau}$ variant as before. Thus for  $L^{12/5}$ estimate  we have
\begin{align*}\| \Big[ \sum_{\tau \in E}(\phi_{\tau}c_{\tau})^2 \Big]^{1/2}\|_{L^{12/5}(B_R)}& \lesssim  \| \Big[ \sum_{\tau \in E}(\phi_{\tau}b_{\tau})^2 \Big]^{1/2}\|_{L^{12/5}(B_R)} \\
& \lesssim 
R^{1/8} \Big[ \sum_{\tau }\int_{B_R}(\phi_{\tau}b_{\tau})^{12/5} \Big]^{5/12}.\end{align*}
Here  on the right hand side $\tau$ ranges over a partition into $R^{-1/2}$ of all of $S$.
 Due to size of our caps we have $b_{\tau} \lesssim R^{-5/2}$. Using this and then conservation of mass  we obtain
\[\lesssim R^{  1/8 +\epsilon}.\]
 We estimate  in $L^{18/7}$ first without using the Kakeya maximal function estimate . First decompose $\phi_{\tau}$ into $\phi_{\tau,\mu}$ exactly as before. Then
\begin{align*}\| \Big[ \sum_{\tau \in E}(\phi_{\tau, \mu}b_{\tau})^2 \Big]^{1/2}\|_{L^{18/7}(B_R)} & \leq R^{1/6} \Big( \sum_{\tau} \|\phi_{\tau,\mu}b_{\tau}\|_{L^{18/7(B_R)}}^{18/7}\Big)^{7/18}\\
&\lesssim R^{1/6}\big(\mu^{-2/21}R^{ \frac{7}{2}+\frac{5}{2}-\frac{5}{2}\cdot \frac{18}{7} +\epsilon }\big)^{7/18}\\
&\lesssim R^{\epsilon}\mu^{-1/27}
\end{align*}

Now we shall estimate using Kakeya maximal function bounds. Again we will first pass to probability measures, and then to the maximal function estimate. Thus we write 
  
\[b_{\tau} \lesssim b_{\tau}\ast R^{-7/2}1_{\tau'}\]
\[(b_{\tau})^2(x) \lesssim R^{-7/2}\int (b_{\tau})^2(y)1_{\tau'}(x-y)dy \]
and estimate using again conservation of mass
\[\int_{B_R}(b_{\tau})^{2}(x) dx \lesssim R^{-3/2}.\]
Hence 
\[(b_{\tau})^2 \lesssim R^{-5}\int 1_{\tau'}(x-y)d\mu_{\tau}\]
So 
\begin{align*}
\| \Big[ \sum_{\tau \in E}(\phi_{\tau, \mu}b_{\tau})^2 \Big]^{1/2}\|_{L^{18/7}(B_R)} & \leq \mu \| \Big[ \sum_{\tau}(b_{\tau})^2 \Big]^{1/2}\|_{L^{18/7}(B_R)}\\ &\lesssim R^{-(\frac{7}{4}+\frac{3}{4})}\mu \|\sum_{\tau}1_{\tau'}(\cdot-y_{\tau})\|_{L^{9/7}(B_R)}^{1/2} \\ &\lesssim R^{-5/2+22/9+\epsilon}\mu.\\
& \lesssim R^{-1/18+\epsilon}\mu.
\end{align*}

Interpolation shows that the improvement for this term is, as in the six-linear case,
\[\frac{2}{735}.\]

Finally comparing the improvements, we get
\[\min(\frac{2}{735},\frac{5}{1764})=\frac{2}{735}\]
and thus we have Theorem 1.


\begin{thebibliography}{99}

\bibitem{bct}
J.Bennett, A. Carbery, T. Tao, On the multilinear restriction and Kakeya conjectures, Acta Math., 196, 261-302  (2006)
\bibitem{b}
J. Bourgain, Besicovitch type maximal operators and applications to Fourier Analysis, Geom.  and  Funct. Anal. 22, 147-187 (1991)
\bibitem{b1}
 J. Bourgain, Harmonic analysis and combinatorics: How much may they contribute to each other?, Mathematics: Frontiers and perspectives, IMU/ Amer. Math. Society, 13-32 (2000)
\bibitem{bg}
J. Bourgain, L. Guth, Bounds on oscillatory integral operators based on multilinear estimates, to appear Geom. Funct. Anal.
\bibitem{f}
C. Fefferman, Inequalities for strongly singular convolution operators, Acta Math. 124, 9-36, (1970) 
\bibitem{la}
 I. Laba, From harmonic analysis to arithmetic combinatorics, Bull. of Amer. Math. Soc., 45,77-115 (2007)  
\bibitem{st}
E.M. Stein, Some problems in harmonic analysis, Harmonic analysis in Eucliden spaces(Proc. Sympos. Pure Math., Williams Coll., Williamstown, Mass., 1978), Part 1, 3-20 (1978)
\bibitem{t}
 T. Tao, A sharp bilinear restriction estimate for paraboloids, Geom. Funct. Anal., 13, 1359-1384 (2003)
\bibitem{t1}
T. Tao, Recent progress on the restriction conjecture, IAS/ Park City Mathematics Series, (2003)
\bibitem{tvv}
T. Tao, A. Vargas, L. Vega, A bilinear approach to the restriction and Kakeya conjectures, J. Amer. Math. Soc., 11, 967-1000 (1998)
\bibitem{tv1}
 T. Tao, A. Vargas, A bilinear approach to cone multipliers I. Restriction Estimates, Geom. Funct. Anal., 10, 185-215 (2000)
 \bibitem{tv2}
 T. Tao, A. Vargas, A bilinear approach to cone multipliers II. Applications, Geom. Funct. Anal., 10, 216-258 (2000)
\bibitem{tomas}
P.Tomas, A restriction theorem for the Fourier transform. Bull. Amer. Math. Soc. 81 477-478 (1975)
\bibitem{w}
T. Wolff,  An Improved bound for Kakeya type maximal functions, Revista Math. Iberoamericana 11, 651-674  (1995)
\bibitem{w1}
T. Wolff, Recent work connected with Kakeya problem, Prospects in Mathematics, 129-162 (1999)
 


\end{thebibliography}
\end{document}